\documentclass[]{interact}

\usepackage{epstopdf}
\usepackage[caption=false]{subfig}
\usepackage[numbers,sort&compress]{natbib}
\bibpunct[, ]{[}{]}{,}{n}{,}{,}
\renewcommand\bibfont{\fontsize{10}{12}\selectfont}
\theoremstyle{plain}
\newtheorem{theorem}{Theorem}[section]

\theoremstyle{definition}
\newtheorem{definition}[theorem]{Definition}

\theoremstyle{remark}

\begin{document}

\title{Spectral Collocation Solutions to Second Order Singular Sturm-Liouville Eigenproblems.}

\author{
\name{C.~I. Gheorghiu\thanks{CONTACT C.~I. Gheorghiu. Email: ghcalin@ictp.acad.ro}} \affil{Tiberiu Popoviciu Institute of Numerical Analysis, Str. Fantanele 57, Cluj-Napoca, Romania}
}

\maketitle
\begin{abstract}
We comparatively use some classical spectral collocation methods as well as highly performing Chebfun algorithms in order to compute the eigenpairs of second order singular Sturm-Liouville problems
with separated self-adjoint boundary conditions. For both the limit-circle non oscillatory and oscillatory cases we pay a particular attention. Some "hard" benchmark problems, for which usual numerical methods (f. d., f. e. m., etc.) fail, are analysed. For the very challenging Bessel eigenproblem we will try to find out the source and the meaning of the singularity in the origin. For a double singular eigenproblem due to Dunford and Schwartz we we try to find out the precise meaning of the notion of continuous spectrum.  For some singular problems only a tandem approach of the two classes of methods produces credible results.
\end{abstract}


\section{Introduction}
There is clearly an interest to develop accurate and efficient methods of solution to singular the Sturm–
Liouville (SL) problems. Our main interest here is to compare the capabilities of Chebfun with those of classical spectral methods 
in solving such problems. The latter
employ basis functions and/or grid points based on Chebyshev and Laguerre polynomials as well as on sinc functions. The effort expended by both classes of methods is also of real interest. It can be assessed in terms of the ease of implementation of the methods as well as in terms of computer resources required to achieve a specified accuracy.

Spectral methods have been shown to provide exponential convergence for a large variety of problems, generally with smooth solutions, and are often preferred.
For details on Chebfun we refer to \cite{DBT, DHT, DHT1, TBD, LNT}  and for Chebyshev collocation (ChC) and Laguerre-Gauss-Radau collocation (LGRC) we refer among other sources to our contributions \cite{Cig14, Cig18} and \cite{WR}. For problems on unbounded domains sinc collocation (SiC) proved to be well suited. Moreover, this method has given excellent results recorded in our contribution \cite{Cig18} and in the works cited there.

We will argue that Chebfun would provide a greater flexibility in solving such problems than the classical spectral methods. This fact is fully true for regular problems. A Chebfun code contains a few lines in which the differential operator is defined along with the boundary conditions and then a subroutine to solve the algebraic eigenproblem. It provides useful information on the optimal order of approximation of eigenvectors and the degree to which the boundary conditions have been satisfied.

Unfortunately, in the presence of singularities, the maximum order of approximation of the unknowns can be reached ($N \geq 4000$) and then Chebfun issues a message that warns about the possible inaccuracy of the results provided.

We came out of this tangle using alternative classical spectral methods. In this way, when we had serious doubts about the accuracy of the solutions given by Chebfun, we managed to establish the correctness of the numerical results.

For two very challenging eigenproblems we have paid a particular attention. The first one is the so called Dunford Schwartz, with two singularities, in origin and at infinity and with a spectrum involving a discrete and a continuous part. The latter is the classical Bessel one for which we have tried to deeply understand the singularity in origin.

A Chebfun code and a MATLAB ChC code are provided in order to exemplify. With minor modifications they could be fairly useful for various numerical experiments.

The structure of these works is as follows. In Section \ref{SSLP} we recall some specific issues for the singular SL problems and in Section \ref{SC} we comment on the Chebfun structure and the classical spectral methods (differentiation matrices, enforcing boundary conditions, etc.). The fourth Section \ref{numerics} is the central section of the paper. Here we analyze no less than eight benchmark problems. In order to separate the "good" from the "bad" eigenvalues we use their relative drift with respect to some parameters. Of all these, the one for which we obtain the most surprising results is the Bessel one. Our analysis is a fairly natural one because the problem itself involves the parameter on which the eigenvalues depend. Generally this is not the case because whenever there is no analytical way to make progress in computing the eigenvalues an artificial parameter is introduced (see \cite{BO} and \cite{LNTex}).  We end up with Section \ref{conclusions} where we underline some conclusions and suggest some open problems.

\section{A class of singular Sturm-Liouville eigenproblems} \label{SSLP}
The Sturm-Liouville problem is to find eigenvalues $\lambda \in \mathbb{C}$ and
eigenfunctions $u(x)$, generally complex valued, satisfying the differential equation:
\begin{equation}
-\left( p\left( x\right) u^{\prime }\right) ^{\prime }+q\left( x\right)
u=\lambda r\left( x\right) u,\ -\infty \leq a<x<b\leq \infty ,  \label{SLeq}
\end{equation}%
where $p,\ p^{\prime },\ q$ and $r$ are continuous on the open interval $\left(
a,b\right) $ and $p\left( x\right) >0,$ $q\left( x\right) >0$ on $\left(
a,b\right).$

With the notations from \cite{PFX} (see also \cite{BEZ}) we observe that this allows either endpoint to be \textit{regular} or \textit{singular}. 

An endpoint $e,$ is regular for (\ref{SLeq}) if and only if
\begin{itemize}
    \item $e$ is finite, and
    \item $1/p(x)$, $q(x)$ and $r(x)$ are absolutely integrable near $x=e$.
\end{itemize}
Otherwise, $x = e$ is called a singular endpoint of equation (\ref{SLeq}).
The boundary conditions must be provided at regular endpoints. At singular points we
assume \textit{Friedrichs boundary conditions} are chosen whenever $\lambda$ is in a range where the
equation is nonoscillatory. This is equivalent to selection of the principal solution at
the nonoscillatory endpoint in both the limit point and limit circle cases. 
 
Some software packages have been designed over time to solve various singular SL problems. The most important would be 
SLEIGN and SLEIGN2, SLEDGE, SL02F and the SLDRIVER interactive package
which supports exploration of a set of Sturm-Liouville problems with the four previously mentioned packages.

In \cite{PF} (see also \cite{PFX}) the authors designed the software package SLEDGE. They  observed that for a class of
singular problems their method either fails or converges very slowly. Essentially, the numerical
method used in this software package replaces the coefficient functions $p(x)$, $q(x)$, and $r(x)$ by step function approximation. Similar behavior has been observed on the NAG code SL02F introduced in \cite{MP} and \cite{PM} as well as on the packages SLEIGN 
and SLEIGN2 introduced in \cite{BEZ} and \cite{BGKZ}.

The main purpose of this paper is to argue that Chebfun, together with the spectral collocation methods, can be a very feasible alternative to these software packages regarding accuracy. In addition these methods can calculate exactly the whole set of eigenvectors and provide some details on the accuracy of the results provided.
 
Towards this end the equation (\ref{SLeq}) can be rewritten in the form%
\[
-\left( p\left( x\right) u^{\prime }\right) ^{\prime }=\tau \left( x\right)
p\left( x\right) u,
\]%
where the quotient $\tau \left( x\right) $ is defined by%
\begin{equation}
\tau \left( x\right) :=\frac{\lambda r\left( x\right) -q\left( x\right) }{%
p\left( x\right) }.  \label{tau}
\end{equation}%
In \cite{PFX} the authors observe that when $x$ approaches an endpoint $e,$
which can be $a$ or $b,$ and the quotient $\tau \left( x\right) $ becomes
positive and unbounded near that endpoint the radial Pr\"{u}fer coordinate
(as a solution of a linear differential equation) will have a rapidly
increasing behavior which can lead to numerical difficulties.

Motivated by this discussion, the authors of \cite{PFX} define the class of hard problems.

\begin{definition}Equation (\ref{SLeq}) is called 'hard" at $\lambda $ near $x = e$ for Sturm-Liouville solvers iff equation (\ref{SLeq}) is nonoscillatory for $\lambda $
at $x = e$ with $e$ finite, and
$\tau \left( x\right) \rightarrow \infty $ as $x\rightarrow e.$
\label{def}
\end{definition}

\section{Spectral collocation (ChC, LGRC, SiC)} \label{SC}

\subsection{Chebfun}
The Chebfun system, in object-oriented MATLAB, contains algorithms which amount to
spectral collocation methods on Chebyshev grids of automatically determined resolution.  Its properties are briefly summarized in \cite{DHT}. In \cite{DBT} the authors explain that chebops are the fundamental Chebfun tools for solving ordinary differential (or integral) equations. One may then use them as
tools for more complicated computations that may be nonlinear and may involve
partial differential equations. This is analogous to the situation in MATLAB itself,
and indeed in computational science generally, where the basic tools are linear
and vector-oriented but they are exploited all the time to solve nonlinear and
multidimensional problems. The implementation of chebops combines the numerical analysis idea of spectral
collocation with the computer science idea of \textit{lazy or delayed evaluation}. The grammar of chebops along with a lot of illustrative  examples is displayed in the above quoted paper as well as in the text \cite{TBD}. Thus it is fairly clear what they can do.

Two eigenproblems, namely a simple Mathieu with periodic boundary conditions, and an Orr-Sommerfeld associated with the eigenvalue instability of plane Poiseuille fluid flow, are also fairly accurate solved.

Actually we want to show in this paper that Chebfun along with chebops can do much more, i.e., can accurately solve highly singular SL eigenproblems. 
 
\subsection{ChC and LGRC}
In all spectral collocation methods designed we have used the differentiation matrices from the seminal paper \cite{WR}. We preferred this MATLAB differentiation suite for the accuracy, efficiency as well as for the ingenious way of introducing various boundary conditions.

In order to impose (enforce) the boundary conditions we have used the \textit{boundary bordering} as well as the \textit{basis recombination}. A very efficient way to accomplish the boundary bordering is available in \cite{JH} and is called \textit{removing technique of independent boundary conditions}. We have used this technique in the large majority of our papers except \cite{GigIsp} where the latter technique has been employed. In the last quoted paper a modified Chebyshev tau method based on basis recombination has been used in order to solve an Orr-Sommerfeld problem with an eigenparameter dependent boundary condition.

In \cite{CIGHPR} we have solved some multiparameter (MEP) eigenproblems which come from separation of variables, in several orthogonal coordinate systems, applied to the Helmholtz, Laplace, or Schrödinger equation. Important cases include Mathieu’s system, Lamé’s system, and a system of spheroidal wave functions.
We show that by combining spectral collocation methods, ChC and LGRC, and new efficient numerical methods for algebraic MEPs, it is possible to solve such problems both very efficiently and accurately. We improve on several previous results available in the literature, and also present a MATLAB toolbox for solving a wide range of problems.

\section{Numerical benchmark problems and discussions} \label{numerics}

\subsection{The Legendre eigenproblem} \label{Leg_polyn}
The Legendre equation reads%
\begin{equation}
-\left( \left( 1-x^{2}\right) u^{\prime }\right) ^{\prime }+\frac{1}{4}%
u=\lambda u,\ -1<x<1.  \label{Leg_eq}
\end{equation}%
In \cite{Pleijel} the author observes that this equation is Weyl's limit-circle type over the
interval $\left( -1,\ 1\right) $. According to Weyl's theory, a symmetric
boundary condition must be added to the differential equation in in order to
define a selfadjoint operator in $L^{2}\left( -1,\ 1\right) $. The symmetric
boundary restriction under which the Legendre polynomials are the eigenfunctions is%
\begin{equation}
\int \left( 1-x^{2}\right) \left\vert u^{\prime }\right\vert ^{2}dx<\infty ,
\label{sym_bc}
\end{equation}%
where the integral is considered over any neighbourhood of the endpoints $x=-1$ and $x=1$. Consequently, we attach to
equation (\ref{Leg_eq}) the Friedrichs boundary conditions%
\begin{equation}
\lim_{x\rightarrow \pm 1}\left[\left( 1-x^{2}\right) u^{\prime }\left( x\right)\right] =0.
\label{Leg_bc}
\end{equation}

On slightly different considerations in \cite{BEZ}, the same boundary
conditions are used in order to find the Legendre polynomials as
eigenfunctions.
\begin{figure}
\centering
\includegraphics[scale=0.75]{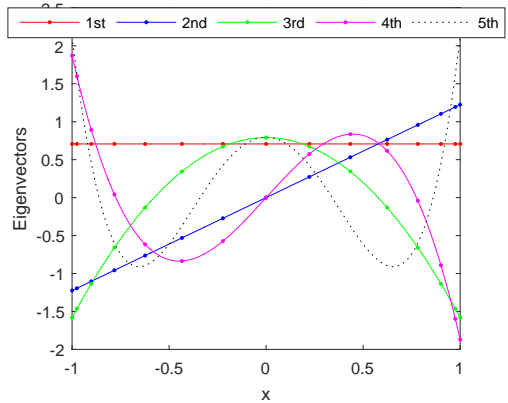}
\caption{\fontsize{9}{9}\selectfont The first five eigenvectors of problem (\ref{Leg_eq})-(\ref{Leg_bc}) at interpolation points computed by Chebfun.}
\label{Leg_pol}
\end{figure}
The boundary value problem (\ref{Leg_eq})-(\ref{Leg_bc}) is the classical
case whose eigenfunctions are the classical Legendre polynomials and whose eigenvalues are known to be:%
\[
\lambda _{n}=n(n+1)+1/4;\ n=0,1,2,3,\ldots \ .
\]%
A simple Chebfun code computes these values with machine precision and also
provides accurate eigenvectors approximating the Legendre polynomials. The first
five of them are depicted in Fig. \ref{Leg_pol}. The code also verifies the restriction (\ref{sym_bc}) which is plainly fulfilled.

It is also important to notice that except for the first eigenvalue $\lambda
_{0}=1/4,$ the above problem could classified as hard in accordance with Definition \ref{def}.
\begin{figure}
\centering
\subfloat[First four eigenvectors of Latzko-Fichera eigenproblem.]{%
\resizebox*{5cm}{!}{\includegraphics{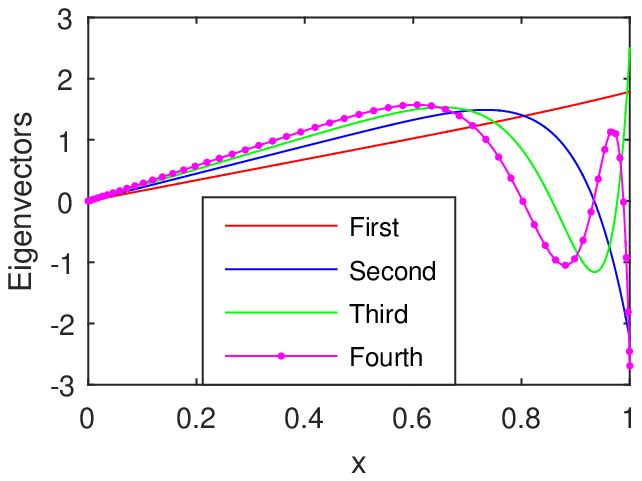}}}\hspace{5pt}
\subfloat[Te Chebyshev coefficients of the eigenvectors.]{%
\resizebox*{5cm}{!}{\includegraphics{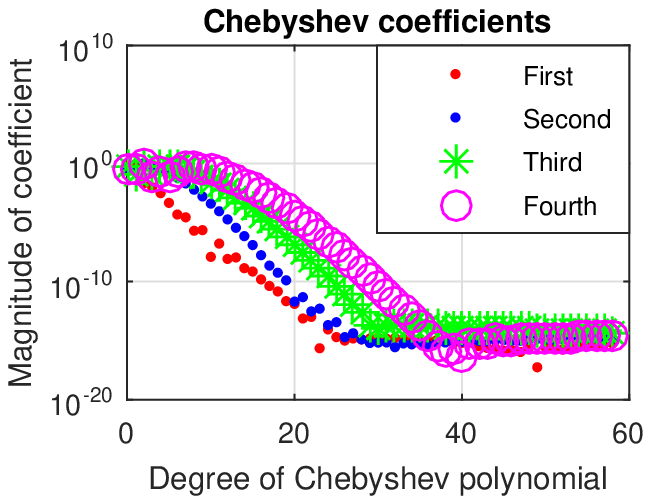}}}
\caption{Chebfun outcomes for Latzko-Fichera eigenproblem. } \label{L_F_fig}
\end{figure}
Moreover, throughout this paper we will say that eigenfunctions will be numerically approximated with eigenvectors and we will draw them.

\subsection{The Latzko-Fichera eigenproblem} \label{L_F_problem}

Latzko (1921) and Fichera (1976) (see \cite{BEZ}) have considered the following singular eigenproblem
\begin{subequations} 
\begin{equation}
-\left( \left( 1-x^{7}\right) u^{\prime }\right) ^{\prime }=\lambda x^{7}u,\
0<x<1,
\label{L_Fichera}
\end{equation}
\begin{equation}
u\left( 0\right) =0,\ \ \lim_{x\rightarrow 1}\left[\left( 1-x^{7}\right) u^{\prime
}\left( x\right)\right] =0.%
\label{L_Fichera_bc}
\end{equation}
\end{subequations}
\begin{table} \centering%
\begin{tabular}{|l|l|}
\hline
$j\:$ & $\:\lambda_{j}$ computed by Chebfun  \\ \hline
$0\:$ & $\:8.727470352650549e+00$ \\ \hline
$1\:$ & $\:1.524230708786303e+02$ \\ \hline
$2\:$ & $\:4.350633321758573e+02$ \\ \hline
$3\:$ & $\:8.556857252681226e+02$ \\ \hline
$4\:$ & $\:1.414142820954995e+03$ \\ \hline
$5\:$ & $\:2.110387972308661e+03$ \\ \hline
\hline
\end{tabular}%
\caption{The first six eigenvalues of problem (\ref{L_Fichera})-(\ref{L_Fichera_bc}).}
 \label{L_Feigs}%
\end{table}%
In the Table \ref{L_Feigs} we report the first six computed eigenvalues.  They are in very good accordance with those reported in \cite{BEZ}. For instance the first eigenvalue coincides up to the sixth digits with that computed by SLEIGN2. In accordance with Definition \ref{def} the problem (\ref{L_Fichera})-(\ref{L_Fichera_bc}) is again hard with respect to all its eigenvalues.

\subsection{A heavy rod-like body}\label{rod}

This problem models a heavy rod-like body with variable cross-section buckling under its own weight. We consider now the singular eigenproblems of the form
\begin{equation}
\left\{
\begin{array}{c}
-\left( A\left( x\right) u^{\prime }\right) ^{\prime }+\left(\gamma-\lambda\right) u\left(
x\right) =0,\ x\in (0,1], \\
\lim_{x\rightarrow 0}A\left( x\right) u^{\prime }\left(
x\right) =0, u\left( 1\right) =0, %
\end{array}%
\right.   \label{P_rod}
\end{equation}%
where $A\in C^{1}\left( \left[ 0,1\right] \right) $ with $A\left( x\right)
>0,$ $x\in (0,1],$ $\int_{0}^{1}A\left( x\right) u^{\prime }\left( x\right)
^{2}dx\,<\infty $ and there is a constant $L\in \left( 0,\infty \right) $
such that $$\lim_{x\rightarrow 0}A\left( x\right) /x^{p}=L,$$ $p\geq 0$ being the\textit{\ tapering} parameter.
\begin{figure}
\centering
\includegraphics[scale=0.98]{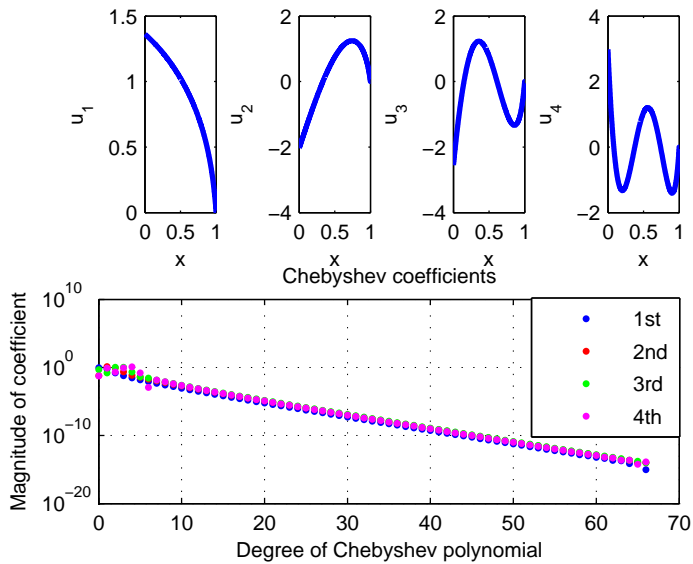}
\caption{\fontsize{9}{9}\selectfont The first four eigenvectors $u_{1},\ldots,u_{4},$ (upper panels) and in a
log-linear plot their Chebyshev coefficients computed by Chebfun (lower panel).}
\label{rod_vect}
\end{figure}
\begin{table} \centering%
\begin{tabular}{|l|l|}
\hline
$j\:$ & $\:\lambda_{j}$ computed by Chebfun  \\ \hline
$0\:$ & $\:1.063402823775151e+00$ \\ \hline
$1\:$ & $\:9.757849576315739e+00$ \\ \hline
$2\:$ & $\:2.575153869531348e+01$ \\ \hline
$3\:$ & $\:4.890432732322349e+01$ \\ \hline
$4\:$ & $\:7.920090041934151e+01$ \\ \hline
$5\:$ & $\:1.166360305002941e+02$ \\ \hline
\hline
\end{tabular}%
\caption{The first six eigenvalues of problem (\ref{P_rod}).}
 \label{P_rod_eigs}%
\end{table}
We perform all numerical experiments using $\gamma:=0$, $A\left(x\right):=log\left(1+sin\left(3x\right)\right)$ and then $p=1$ and $L=3.$ The first four computed eigenvectors are displayed in Fig. \ref{rod_vect}. With respect to the accuracy of our computations, from the lower panel of Fig. \ref{rod_vect}, we observe that the series solution for each eigenvector converges rapidly and smoothly. By the time we reach an approximation order of about $65$, the accuracy is around $15$ digits, and the computed Chebyshev series is truncated there. Thus the attained accuracy is of machine precision for each of the first four vectors.

The first six eigenvalues computed by Chebfun are displayed in Table \ref{P_rod_eigs}. 
For $\gamma:=0$, in \cite{S} the author provides a rigorous spectral theory of the unbounded linear operator involved.
The main conclusion is that the
spectral properties of the problem for tapering of order $p = 2$ are very different from
what occurs for $p < 2$. For $p = 2$, there is a non-trivial essential spectrum and
possibly no eigenvalues, whereas for $p < 2$, the whole spectrum consists of a sequence
of simple eigenvalues (see also \cite{Castro}). We plainly confirm numerically this latter statement.

At the end of this subsection it is important to note that corresponding to each eigenvalue the problem is hard, and the case $\gamma \neq 0$ does not involve supplementary computing complications.

Moreover, it is worth mentioning that until now we have not been able to introduce Friedrichs type boundary conditions when we have applied the classical spectral methods. In contrast, the introduction of these conditions with Chebfun is extremely simple and direct.

\subsection{A Sturm-Liouville eigenvalue problem
with an interior singularity (Boyd eigenproblem).}\label{Boyd}

In \cite{EGZ} the authors consider the following singular, but slightly less general eigenproblem than (\ref{SLeq}), namely
\begin{equation}
-u^{\prime \prime }+q\left( x\right) u=\lambda u,\ -\infty <a<x<b<\infty ,
\label{SLeq1}
\end{equation}%
with $q$ real-valued. It is called singular if $q\in L\left[ \alpha ,\ \beta %
\right] \ $\ for any $\left[ \alpha ,\ \beta \right] \subset \left[ a,\ b%
\right] $ and $q\notin L\left[ a,\ b\right] .$ 
For two real numbers $a$ and $b$, such that $-\infty
<a<0<b<\infty ,\ q\left( x\right) :=-1/x,$ and homogeneous Dirichlet boundary conditions 
\begin{equation}
u\left( a\right) =u\left( b\right) =0,\label{Hom_D}
\end{equation}
attached to (\ref{SLeq1}), in \cite{EGZ} the authors show that this eigenproblem
has a unique solution which has a discrete real spectrum $\left\{ \lambda
_{n},\ n\in\mathbb{N}_{0}\right\} $ where%
\[
-\infty <\lambda _{0}\leq \lambda _{1}\leq \lambda _{2}\leq \ldots \leq
\lambda _{n}\leq \lambda _{n+1}\leq ...,\ \lim_{n\rightarrow \infty }\lambda
_{n}=\infty.
\]%
The corresponding real-valued eigenfunctions $\left\{ u_{n},\ n\in \mathbb{N}
_{0}\right\} $ satisfy the equation (\ref{SLeq1}) and the Dirichlet
boundary conditions (\ref{Hom_D}) and the following properties for $n\in 
\mathbb{N}_{0}$:
\begin{itemize}
\item $u_{n}:\left[ a,\ b\right] \rightarrow 
\mathbb{R}
$ and $u_{n}\in C\left[ a,\ b\right] ;$

\item $u_{n}\in C^{2}[a,\ 0)\cup C^{2}(0,b];$

\item $u_{n}(0)=0;$

\item in general $u_{n}^{\prime }(0^{-})\neq $ $u_{n}^{\prime }(0^{+})$ but both limits are finite;

\item $\int_{a}^{b}u_{n}\left( x\right) u_{m}\left( x\right)dx =\delta _{n,m},$
for $n\neq m\in \mathbb{N}_{0};$
\item $\left\{ u_{n},\ n\in \mathbb{N}
_{0}\right\} $ is a complete orthonormal set in the Hilbert space $L^{2}%
\left[ a,\ b\right] ;$
\item degeneracy may occur but no eigenvalue can have multiplicity greater than two.
\end{itemize}
The Chebfun code has failed to solve this problem in form (\ref{SLeq1})-(\ref{Hom_D}). However, in order to solve it numerically we use 
a perturbed version of $q\left(x\right)$, namely
\begin{equation}
q\left( x,\varepsilon \right) :=\frac{x}{x^{2}+\varepsilon ^{2}},
\label{Q_reg}
\end{equation}%
in the sense of distributions, i.e., $\lim_{\varepsilon \rightarrow
0}q\left( x,\varepsilon \right) =\delta \left( \frac{1}{x}\right) .$ With this approximation of the original coefficient $q\left(x\right)$ the eigenproblem becomes a \textit{regular} one. It is
temping to hope that the the perturbed set of eigenvalues $\left\{ \lambda
_{n}\left( \varepsilon \right) ,\ n\in \mathbb{N}
_{0}\right\} $ and of eigenvectors $\left\{ u_{n}\left( x,\varepsilon
\right) ,\ n\in \mathbb{N}
_{0}\right\} $ approaches in some sense the corresponding sets of the
original problem when $\varepsilon \rightarrow 0.$

Actually with Chebfun computations we have used $\varepsilon:=1.e-06$ in (\ref{Q_reg}) on the integration interval $\left[-10, 10\right]$.

\subsubsection{Sinc collocation results}\label{SiC}

\begin{figure}
\centering
\subfloat[First four eigenvectors of Boyd eigenproblem computed by SiC.]{%
\resizebox*{5cm}{!}{\includegraphics{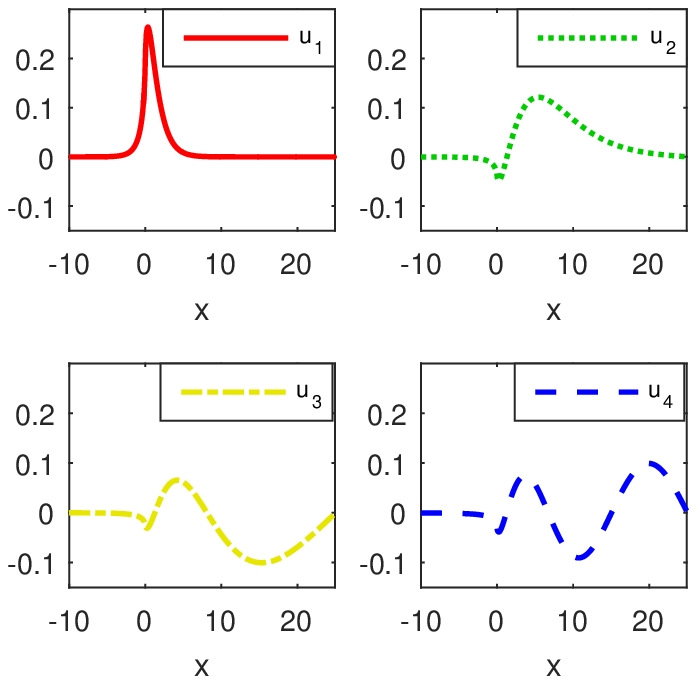}}}\hspace{5pt}
\subfloat[The sinc coefficients of the eigenvectors.]{%
\resizebox*{5cm}{!}{\includegraphics{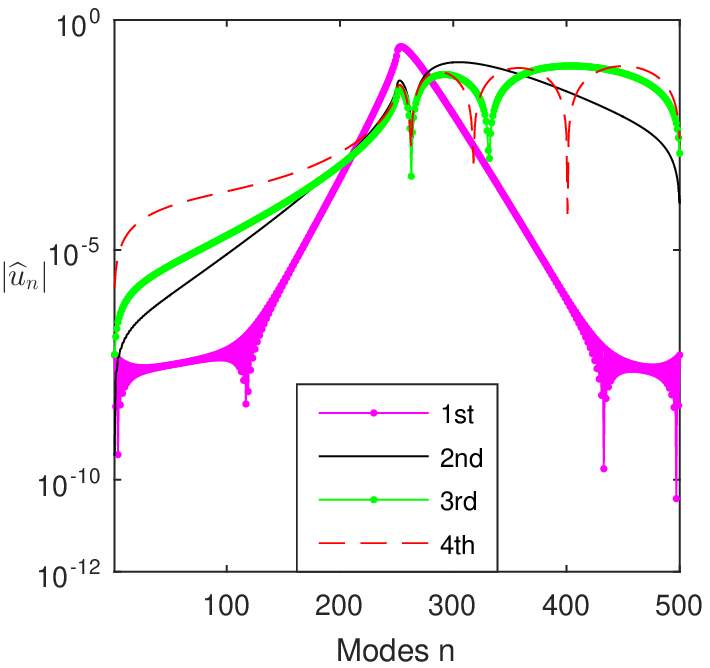}}}
\caption{SiC outcomes for Boyd eigenproblem. } \label{Boyd_fig}
\end{figure}

The unsatisfactory results obtained with Chebfun have prompted us to use an alternative method.
In order to implement the SiC we have used the order of approximation $N:=500$ and the scaling factor $h:=0.1$. To the parameter $\varepsilon$ in (\ref{Q_reg}) has been assigned the value $1.e-06$. The first four vectors are displayed in Fig. \ref{Boyd_fig} (a). They are continuous in origin but with discontinuous derivatives. In Fig. \ref{Boyd_fig} (b) we show the behaviour of their coefficients in SiC formulation. For the first eigenvector the slop looks smooth, even linear. A rouding-off plateau is observed under $10^{-6}$. For the following eigenvectors things get worse gradually.
\begin{table} \centering%
\begin{tabular}{|l|l|l|l|}
\hline
$j\:$ &\:$\lambda_{j}$ computed by SiC & \:$\lambda_{j}$ computed by Chebfun &\:$ \lambda_{j}$ computed in \cite{EGZ}\\ \hline
$0\:$ &$\:-9.606833569044633e-01 $& $\:-9.794200447802075e-01$ & $\: -0.9841$\\ \hline
$1\:$ &$\:-1.095983388802928e-01 $& $\:-7.751928355451891e-02$ & $\: -0.0778$\\ \hline
$2\:$ &$\:-2.965096161409372e-02 $& $\: 2.732439098153830e-01$ & $\:  0.2727$\\ \hline
$3\:$ &$\: 4.851960256747148e-02 $& $\: 3.092576560488862e-01$ & $\:  0.3092$\\ \hline
$4\:$ &$\: 9.064724518382954e-02 $& $\: 6.754139419849143e-01$ & $\:  0.6754$\\ \hline
$5\:$ &$\: 1.630157312824040e-01 $& $\: 8.403112262426201e-01$ & $\:  0.8396$\\ \hline
\hline
\end{tabular}%
\caption{The first six eigenvalues of Boyd eigenproblem computed by three different methods.}
 \label{Boyd_eigs}%
\end{table} 
In Table \ref{Boyd_eigs} we report the first six eigenvalues computed by SiC, Chebfun and by the SLEIGN package according to \cite{EGZ}. In the latter case the authors worked with $\varepsilon:=1.e-03$ on the interval $\left[-10, 10\right]$. There is a good agreement between the results obtained with Chebfun and SLEIGN.

\subsection{Fokker-Planck eigenproblem} \label{F_P_eig}

In \cite{BDS} the author considers a problem of type (\ref{SLeq1}) with 
\begin{equation}
q\left( x\right) :=\frac{x^{6}}{4}-\frac{3x^{2}}{2},\ -\infty <x<\infty 
\label{q_F_P}
\end{equation}%

Symmetrically to the origin, at large distances, we will ask that the
solutions be bounded, i.e., $u\left( \pm l\right) :=0$ with $l \gg 0.$ Along
with this interval truncation hypothesis we have used Chebfun and then SiC. 

The eigenvalues obtained using Chebfun are 
compared with those obtained in \cite{BDS} when the interval of integration has length $l:=4.$ They are reported in Table 
\ref{F_P_eigs} and are in an excellent agreement. We have to mention that the ground state eigenvalue for Fokker-Planck eigenproblem is $\lambda_{0}=0$. This fact is confirmed by our computation as we have got $\lambda_{0}=O(10^{-14})$.  
\begin{table} \centering%
\begin{tabular}{|l|l|l|}
\hline
$j\:$ &\:$\lambda_{j}$ computed by Chebfun & \:$\lambda_{j}$ computed in \cite{BDS} \\ \hline
$1\:$ &$\:1.368592520979542e+00 $& $\:1.36860$ \\ \hline
$2\:$ &$\:4.453709163213802e+00 $& $\:4.45371$ \\ \hline
$4\:$ &$\:1.275806953296428e+01 $& $\: 12.7581$ \\ \hline
$6\:$ &$\: 2.349440842267923e+01 $& $\: 23.4944$ \\ \hline
$10\:$ &$\: 5.061402223182223e+01 $& $\: 50.6140$ \\ \hline
$20\:$ &$\: 1.432321465884990e+02$& $\: 143.232$ \\ \hline
$30\:$ &$\: 2.631594491758098e+02 $& $\: 263.159 $ \\ \hline
\hline
\end{tabular}%
\caption{Some eigenvalues of Fokker-Planck eigenproblem computed by Chebfun and by a pseudospectral method
based on nonclassical polynomials.}
 \label{F_P_eigs}%
\end{table} 

From Fig. \ref{F_P_eigs} it is clear that Chebfun has used an approximation of order $80$ in order to compute this eigenvector. 
Our numerical experiments have shown that slightly increased values of $N$ were used in order to calculate the other higher order vectors.
\begin{figure}
\centering
\subfloat[The fourth eigenvector of Fokker-Planck eigenproblem.]{%
\resizebox*{5cm}{!}{\includegraphics{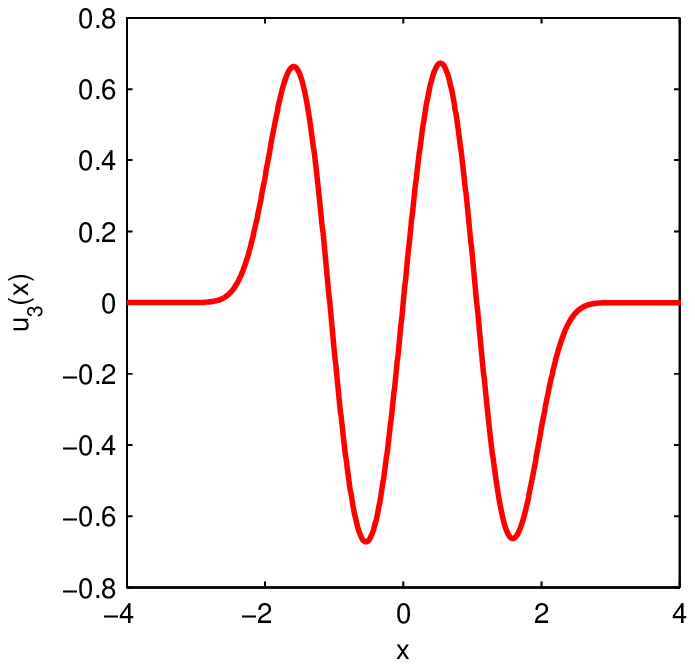}}}\hspace{5pt}
\subfloat[The even and odd Chebyshev coefficients of the fourth eigenvector.]{%
\resizebox*{5cm}{!}{\includegraphics{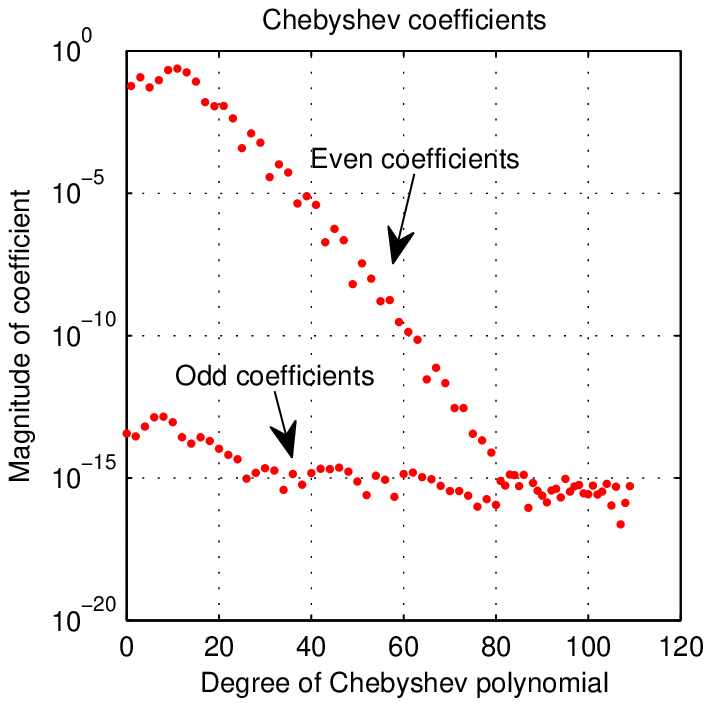}}}
\caption{The fourth eigenvector of Fokker-Planck eigenproblem along with its Chebyshev coefficients. } \label{F_P_fig}
\end{figure}

In order to trust the computed eigenvalues we use the 
concept of \textit{relative drift with respect to a parameter} say $\alpha$. This important concept has been introduced in \cite{Boyd}. An obvious way to achieve the separation of "good" and "bed" eigenvalues is to compute them for for different orders on approximation 
of parameter $\alpha$. Thus the relative drift of the $j th$ eigenvalue is defined as (see \cite{Boyd})

\begin{equation}
\delta _{j,relative,\alpha}:=\frac{\left\vert \lambda _{j}^{\alpha_{1}}-\lambda
_{j}^{\alpha_{2}}\right\vert }{\left\vert \lambda _{j}^{\alpha_{1}}\right\vert },\
\alpha_{1}\neq \alpha_{2},\   \label{eig_drift}
\end{equation}
where $\lambda _{j}^{\left( \alpha\right) }$ is the $j th$ eigenvalue, after the
eigenvalues have been sorted, as computed using a  specified value of parameter $\alpha$. Only those whose relative difference or 
\textit{resolution-difference drift} is small, relative to a desired approximation, can be believed. 

\begin{figure}
\centering
\includegraphics[scale=0.65]{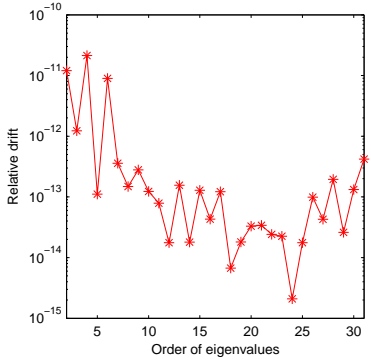}
\caption{\fontsize{9}{9}\selectfont The drift of the eigenvalues $\lambda_{2}$,...,$\lambda_{31}$ of the Fokker-Planck eigenproblem.}
\label{F_P_drift}
\end{figure}

The drift of the eigenvalues $\lambda_{2}$,...,$\lambda_{31}$ of the Fokker-Planck eigenproblem with respect to the length 
$l$ of the interval of integration, i.e., $\alpha:=l$, is computed using formula (\ref{eig_drift}) and is depicted in Fig. \ref{F_P_drift} in a log-linear plot. We have used 
$l:=4$ and $l:=10$ and can conclude that the eigenvalues are correctly computed within an approximation of order $10^{-11}$.

\subsection{Dunford and Schwartz eigenproblem (1963)} \label{D_S_problem}

The generalized hypergeometric equation reads
\begin{equation}
-u^{\prime \prime }+\frac{-242\cosh (x)+241}{4\sinh ^{2}\left( x\right) }%
u=\lambda u,\ 0<x<\infty.  \label{D_S_eq}
\end{equation}%
In \cite{Pr} the author attaches to this equation the boundary conditions%
\begin{equation}
\lim_{x\rightarrow 0}\left[u\left( x\right) +xu^{\prime }\left( x\right)
\right] =0,\ u\left( \infty \right) =0.  \label{D_S_bc}
\end{equation}

The singularity in origin $\left( \sinh ^{2}\left( x\right) =0\left(x^{2}\right) \ \text{as\ }x\rightarrow 0\right) $ leads to the above Friedrichs type boundary condition in origin. Actually, in order to resolve the singularity at $\infty$ we solve this problem on the truncated interval $\left[\epsilon, X\right]$. The best results obtained have been for $\epsilon:=1e-08$ and $X:=15$. In order to improve the convergence of Chebfun we use the option \texttt{'splitting' 'on'.}

The first ten eigenvalues obtained are reported in Table \ref{D_S_eigs}. The first five of them reasonable approximate the exact eigenvalues $\lambda_{n}=-\left(5-n\right)^2$, $n=0,1,2,3,4$. The next five (the second column in Table \ref {D_S_eigs}) suggest the \textit{continuous spectrum} on $\left(0, \infty\right).$
\begin{table} \centering%
\begin{tabular}{|l|l||l|l|}
\hline
$j\:$ & \:Discrete $\lambda_{j}$ by Chebfun & $j\:$ & \:"Continuous spectrum" by Chebfun \\ \hline
$0\:$ & $\:-2.493732084634685e+01$ & $5\:$ & $\:1.578606493636769e-02$\\ \hline
$1\:$ & $\:-1.572675883354506e+01$ & $6\:$ & $\:2.297844931237520e-01$\\ \hline
$2\:$ & $\:-8.878700805769942e+00$ & $7\:$ & $\:6.369715876528872e-01$\\ \hline
$3\:$ & $\:-3.934516542596402e+00$ & $8\:$ & $\:1.133491876230128e+00$\\ \hline
$4\:$ & $\:-8.976614614351384e-01$ & $9\:$ & $\:1.761853804586405e+00$\\ \hline
\hline
\end{tabular}%
\caption{The first ten eigenvalues of Dunford Schwartz eigenproblem (\ref{D_S_eq})-(\ref{D_S_bc}).}
 \label{D_S_eigs}%
\end{table}
Starting with the sixth eigenvalue our numerical experiments really suggest the existence of a \textit{continuous spectrum} on 
$\left(0, \infty \right).$
\begin{figure}
\centering
\subfloat[A zoom in the first four eigenvectors of Dunford Schwartz eigenproblem (\ref{D_S_eq})-(\ref{D_S_bc}).]{%
\resizebox*{5.5cm}{!}{\includegraphics{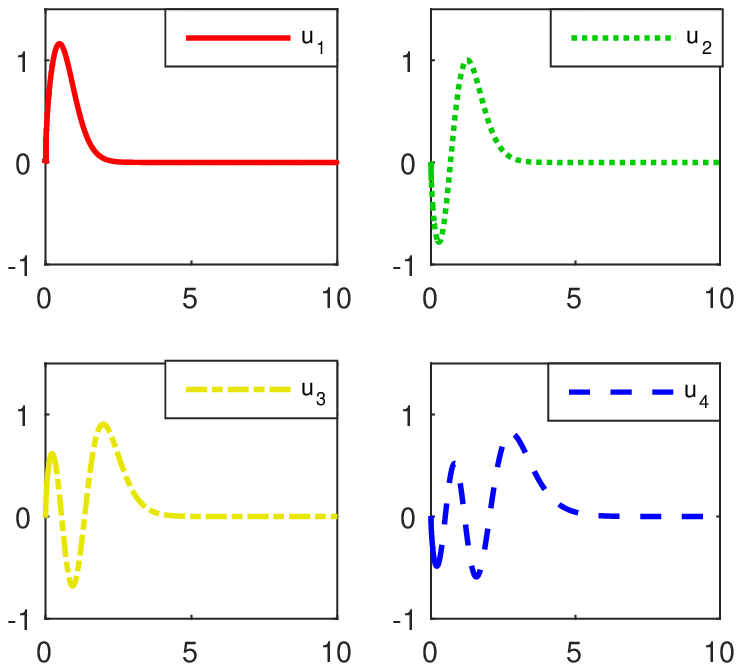}}}\hspace{5pt}
\subfloat[The Chebyshev coefficients of the first eigenvector.]{%
\resizebox*{5.5cm}{!}{\includegraphics{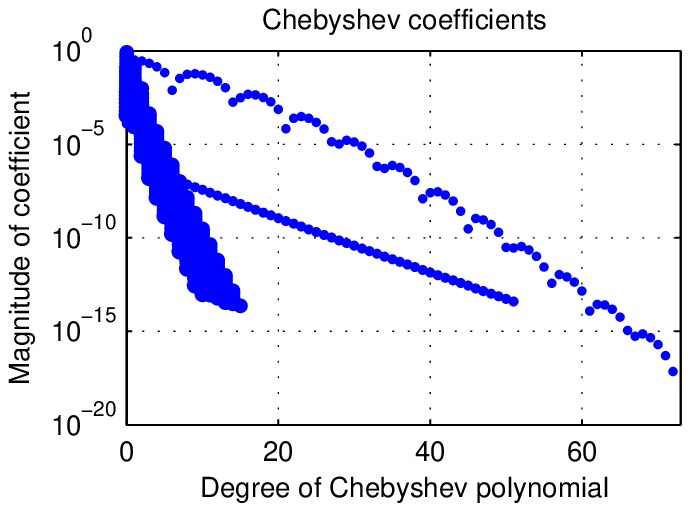}}}
\caption{Chebfun outcomes with respect to eigenvectors for Dunford and Schwartz eigenproblem (\ref{D_S_eq})-(\ref{D_S_bc}). } \label{D_S}
\end{figure}
The first four eigenvectors of Dunford and Schwartz eigenproblem are displayed in Fig. \ref{D_S} a). They look plausible viz. the number of their roots. The behavior of the Chebyshev coefficients of the first eigenvector is shown in Fig. \ref{D_S} b). The odd coefficients drawn below the even ones have a somewhat strange and yet unexplained behavior.

With its two singularities this problem has proved to be one of the most difficult and expensive in terms of computing elapsed time.

\subsection{A nasty $q(x)$ eigenproblem.}\label{nasty_y}

Let's consider now the eigenproblem (\ref{SLeq1})-(\ref{Hom_D}) with 
\begin{equation}
q\left(x\right):=\ln{x}. \label{q_nasty}
\end{equation}
The problem is solved in \cite{Pr} on the interval $\left[0, 4\right]$, as problem (11) in the standards problem set.

\subsubsection{ChC method} 

The problem is of course singular but the classical ChC works extremely fast and accurate. The code consists in finding out of eigenpairs of the matrix $A$ with the simple MATLAB sequence:
\begin{verbatim}
    N=512;                     % The approximation order
    [x,D]=chebdif(N,2);        % Differentiation matrices [WR]
    k=2:N;                     % Kept modes
    xbc=x(k); Deig=D(k,k,2);   % Kept nodes and matrix
    A=-Deig/4+diag(log((xbc+1)*2)); % [0,4] shifted to [-1,1]
    [U,S]=eig(A); [t,o]=sort(diag(S)); S=S(o); U=U(:,o);
    disp(S(1:24))              % Compute,sort and display 
\end{verbatim}
We enforce the boundary condition by \textit{removing technique of independent boundary conditions} introduced in \cite{JH}. Thus, we delete the first and the last rows and columns in the second order Chebyshev differentiation matrix $D(:,:,2)$ and obtain the matrix $Deig$. The nodes $x_{1}$ and $x_{N}$ of the Chebyshev-Gauss-Lobatto system are eliminated as \textit{slaves} and after the problem is solved they are given-back. Actually we simply add two zeros, one in the first and another in the last position of each eigenvector.
\begin{figure}
\centering
\subfloat[The Chebyshev coefficients of the first four eigenvectors computed by ChC.]{%
\resizebox*{6.0cm}{!}{\includegraphics{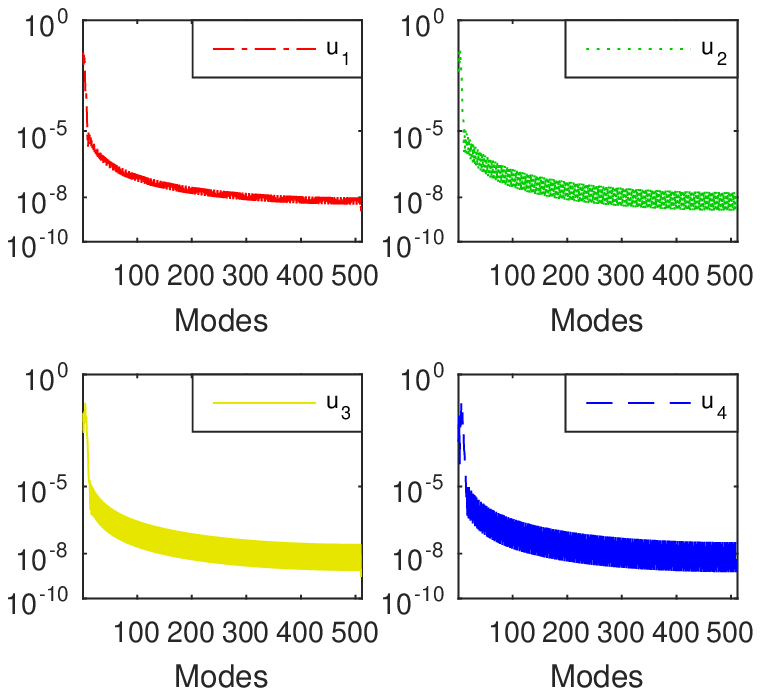}}}\hspace{5pt}
\subfloat[The relative drift with respect to the order of approximation of the first twenty-four eigenvalues, computed with formula (\ref{eig_drift}). We have used $N_{1}:=512$ and $N_{2}:=256$.]{%
\resizebox*{5.0cm}{!}{\includegraphics{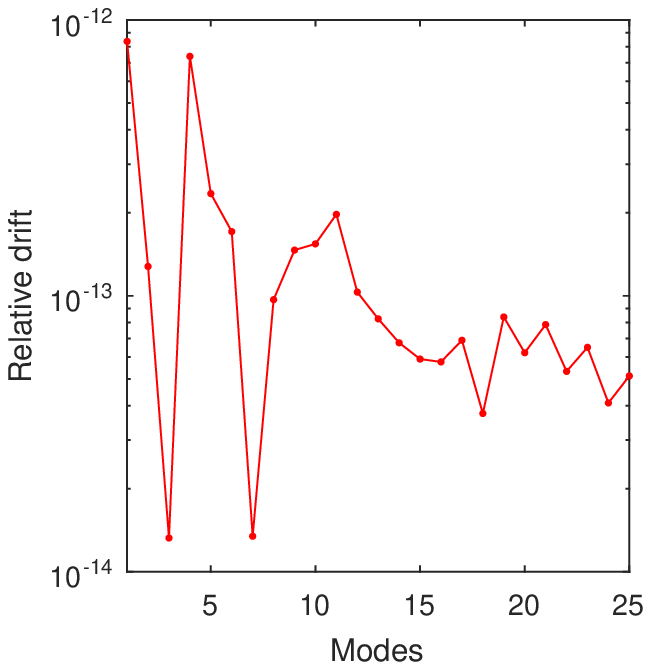}}}
\caption{ChC outcomes with respect to eigenvectors for nasty $q(x)$ eigenproblem \ref{q_nasty}. } \label{ChC_fig}
\end{figure}

\subsubsection{Chebfun method}
\begin{figure}
\centering
\subfloat[First four eigenvectors of nasty $q(x)$ eigenproblem computed by Chebfun.]{%
\resizebox*{5cm}{!}{\includegraphics{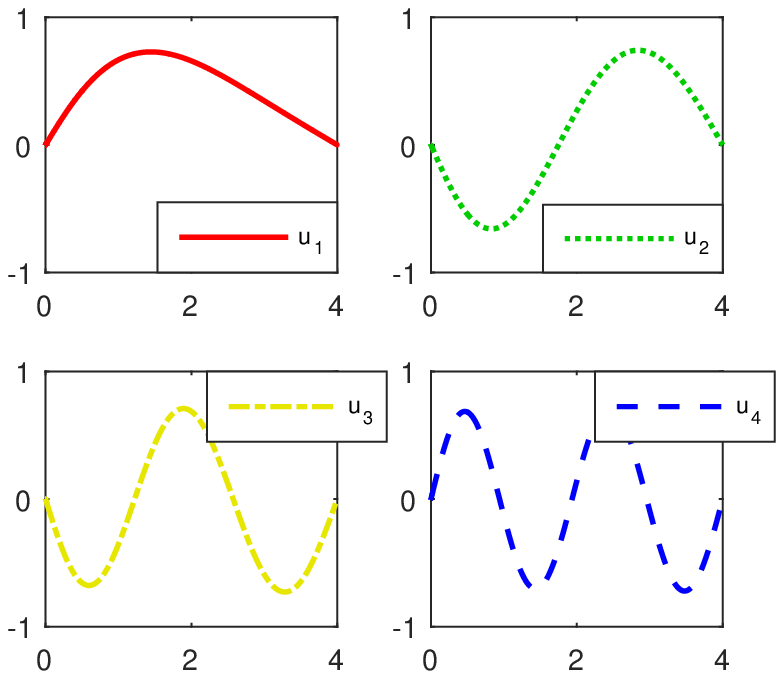}}}\hspace{5pt}
\subfloat[The Chebyshev coefficients of the first four eigenvectors.]{%
\resizebox*{5cm}{!}{\includegraphics{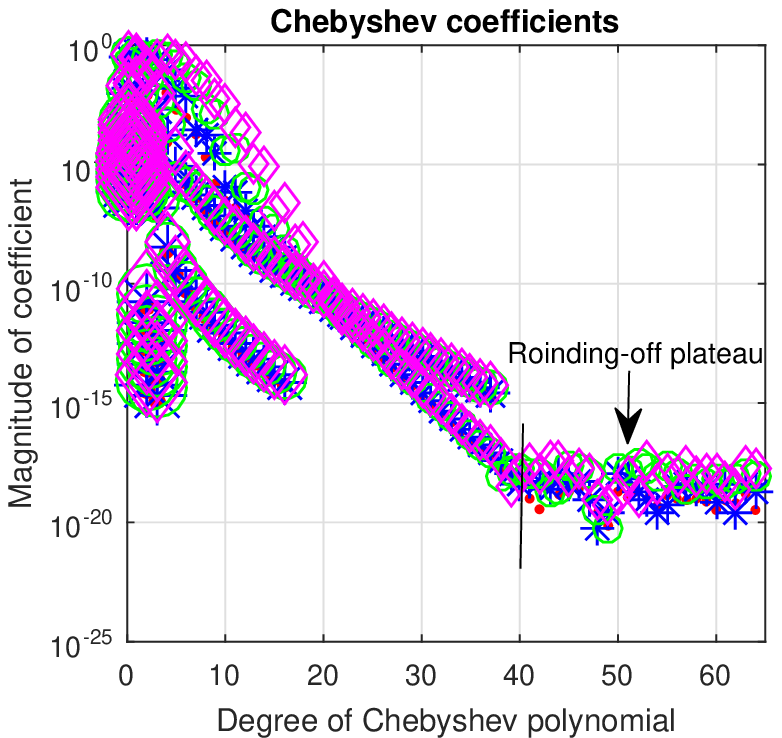}}}
\caption{The Chebfun outcomes for the nasty $q(x)$ eigenproblem (\ref{q_nasty}). } \label{ln_fig}
\end{figure}

In order to avoid the singularity in origin we solve the same problem by Chebfun on the modified (truncated) domain $\left[1.e-08, 4\right]$. The first four eigenvectors are displayed in Fig. \ref{ln_fig} (a). They look fairly smooth and they do not seem affected by the singularity at all. However, their coefficients do not decrease smoothly and this is due to the singularity of $q (x)$. Actually the bandwidth on which they oscillate increases with $n$. This fact is visible in Fig. \ref{ln_fig} (b).
Anyway, no warning has been displayed by Chebfun when the option \texttt{'splitting', 'on'} has been invoked.
\begin{table} \centering%
\begin{tabular}{|l|l|l|l|}
\hline
$j\:$ &\:$\lambda_{j}$ computed by ChC & \:$\lambda_{j}$ computed by Chebfun &\:$ \lambda_{j}$ computed in \cite{Pr}\\ \hline
$0\:$ &$\:1.124816809695236e+00 $& $\:1.124816818756614e+00 $ & $\: 1.1248168097$\\ \hline
$23\:$ &$\:3.557030079371902e+02$& $\:3.557030097207584e+02$ & $\: 385.92821596$\\ \hline
\hline
\end{tabular}%
\caption{The first and the twenty-fourth eigenvalues of nasty eigenproblem (\ref{q_nasty}) computed by three different methods.}
 \label{nasty_eigs}%
\end{table}

\subsection{ChC vs. Chebfun}\label{ChC_Chebfun}

The first and the twenty-fourth eigenvalues of this eigenproblem computed by ChC, Chebfun and reported in \cite{Pr} are displayed in Table 
\ref{nasty_eigs}. If, with respect to the first eigenvalue, we observe that they coincide up to the fifth decimal, this does not happen for the value of twenty-fourth eigenvalue. Thus, in order to separate the "good" eigenvalues from the "bad" ones, i.e., inaccurate eigenvalues, we use again the concept of relative drift, now with respect to the order of approximation $N$, i.e., $\alpha:=N$ in formula (\ref{eig_drift}). Actually, the relative drift of eigenvalues computed by ChC reported in Fig. \ref{ChC_fig} (b) shows a coincidence of at least order $10^{-13}$ for the twenty-fourth eigenvalue computed by ChC. In other words, it means that the relative drift is resolution independent with respect to $N$ for the outcomes of ChC method. Moreover, the method is \textit{numerically stable}. Inspecting Table \ref{nasty_eigs} we can conclude that the eigenvalues computed by ChC and Chebfun are more accurate than those reported in \cite{Pr} especially for large values of $n$.

The eigenvectors computed by the two methods look practically identical, so we only present those obtained through Chebfun in Fig. \ref{ln_fig} a). 
Their Chebyshev coefficients are represented in Fig. \ref{ln_fig} b). They piecewise decrease on a rather regular way up to $N=40$ after which a rounding-off plate is installed. 

The Chebyshev coefficients of the eigenvectors computed by ChC are obtained using the MATLAB code \texttt{fcgltran} from \cite{GvW}. For the first four eigenvectors they are displayed in Fig. \ref{ChC_fig} a). They decrease smoothly in bands (strips) with increasing bandwidth as $n$ increases.

However, comparing the Fig. \ref{ChC_fig} a) and Fig. \ref{ln_fig} b) it becomes very clear that we get a better accuracy with Chebfun with a much smaller effort, with respect to $N,$  but with the price of truncating the domain.

\subsection{The Bessel equation.}\label{Bessel}

This is the equation
\begin{equation}
-u^{\prime \prime }+\frac{c}{x^{2}}u=\lambda u\left( x\right) ,\ 0<x\leq 1,\
c\in \mathbb{R}.  \label{B_eq}
\end{equation}%
The endpoint $1$ is \textit{regular} and $0$ is a \textit{singular} endpoint
for all $c\neq 0.$ Let $c:=\nu ^{2}-1/4.$ In \cite{BEZ} the authors solve two
boundary value problems attached to (\ref{B_eq}). 

For the first one the boundary conditions read
\begin{equation}
\lim_{x\rightarrow 0}\left[ \left( \nu +\frac{1}{2}\right) u\left( x\right)
+xu^{\prime }\left( x\right) \right] =0,\ u\left( 1\right) =0,  \label{NOC}
\end{equation}%
i.e., we are in the the nonoscillatory case $c\neq 0,\ -\frac{1}{4}\leq c<\frac{3}{4}$. The transcendental equation for the BVP (\ref{B_eq})-(\ref{NOC}) is simply 
\[
J_{\nu }\left( s\right) =0,
\]%
where $J_{\nu }$ is the Bessel function of order $\nu.$

For the second one, in the same case, the boundary conditions read%
\begin{equation}
\lim_{x\rightarrow 0}\left[ \left( 1+\frac{1}{2}\ln{x}\right) u\left( x\right)
-\left( x\ln{x}\right) u^{\prime }\left( x\right) \right] =0,\ u\left( 1\right)
=0.  \label{NOC1}
\end{equation}
The transcendental equation for the BVP (\ref{B_eq})-(\ref{NOC1}) is similarly 
\[
J_{-\nu }\left( s\right) =0.
\]
\begin{figure}
\centering
\subfloat[First three eigenvectors of the Bessel eigenproblem computed by Chebfun when $\nu:=0$.]{%
\resizebox*{5cm}{!}{\includegraphics{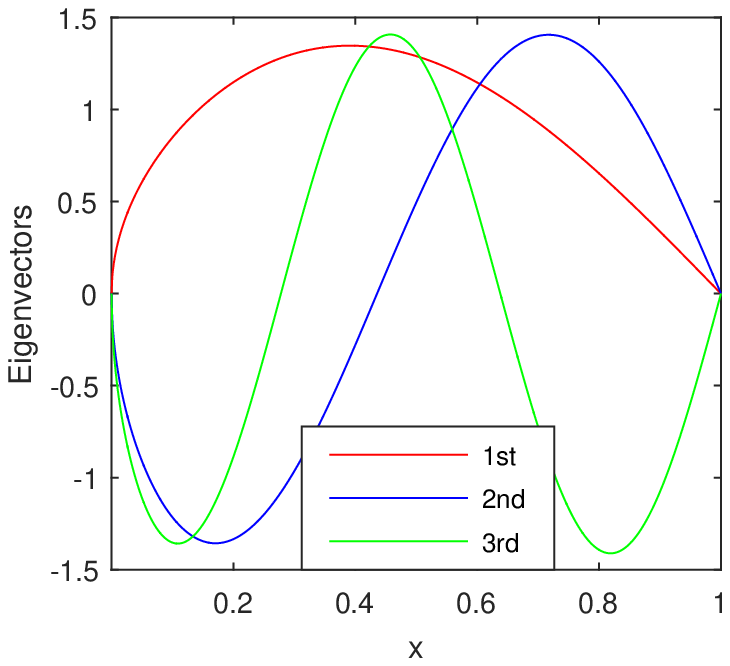}}}\hspace{5pt}
\subfloat[The Chebyshev coefficients of the first three eigenvectors.]{%
\resizebox*{5cm}{!}{\includegraphics{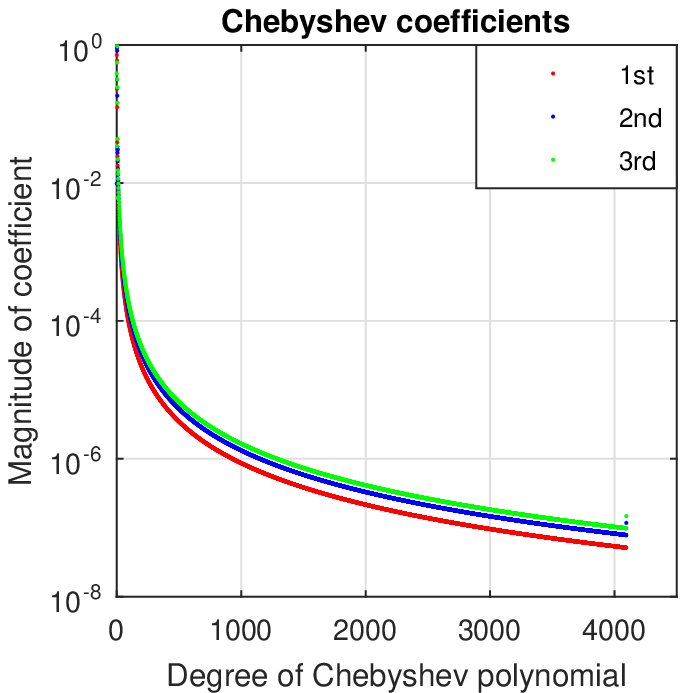}}}
\caption{The Chebfun outcomes for the eigenproblem (\ref{B_eq})-(\ref{NOC}) when $\nu:=0$. } \label{B_0__fig}
\end{figure}
To keep the length of the paper to a reasonable size, we will only report the results of solving the first problem.
As in the vast majority of the problems considered in this paper, we did not find in the literature eigenvectors with which to compare those shown in Fig. \ref{B_0__fig} (a). Regarding their Chebyshev coefficients, we observe from Fig. \ref{B_0__fig} (b) that they descend very smoothly to something of the order $10^{-7}$ that is well above the machine precision. Moreover, Chebfun uses an order of approximation extremely large, i.e, $N>4000$. All these obviously are due to the singularity in the origin.
Listed below are the first three eigenvalues calculated with Chebfun. They coincide very well with the corresponding values frequently cited in the literature. 
\[
\lambda_{0}=2.404996056333427e+00;\:\lambda_{1}=5.520251599508041e+00;\:\lambda_{8}=2.749367004065968e+01.
\]
We also want to observe that the problem consumes a lot of \texttt{cpu} time in the sense that the elapsed time is up to 40-50 times greater than in the case of the Legendre problem, for example (see Subsection \ref{Leg_polyn}).

In \cite{HWV} the author considers a very interesting generalization of problem (\ref{B_eq}). We will reconsider this from a numerical point of view in order to find a deeper understanding of the singularity in the origin.

This new problem reads
\begin{equation}
\left\{ 
\begin{array}{c}
u^{\prime \prime }+\left[ \lambda +\frac{\left( \frac{1}{4}-\nu ^{2}\right) 
}{\left( x-\tau \right) ^{2}}\right] u=0,\ -1<x<1,\ \nu ,\ \tau \in \mathbb{R}, \\ 
u\left( \pm 1\right) =0.%
\end{array}%
\right.   \label{B_HWV}
\end{equation}
For given $\tau >1$ (or $\tau <-1$) we have a regular Sturm-Liouville problem with eigenvalues
\[
0<\lambda _{1}\left( \tau \right) <\lambda _{2}\left( \tau \right) <\ldots .
\]%
Chebfun confirms this statement and solve this regular problem in no time.
It is easy to show that the functions $\lambda _{n}\left( \tau \right) $ are
analytic for the above values of $\tau $ and the author of \cite{HWV} asks
for properties of the analytic continuation of $\lambda _{n}\left( \tau\right) $ into the complex plane.

\begin{figure}
\centering
\includegraphics[width=3.5 in]{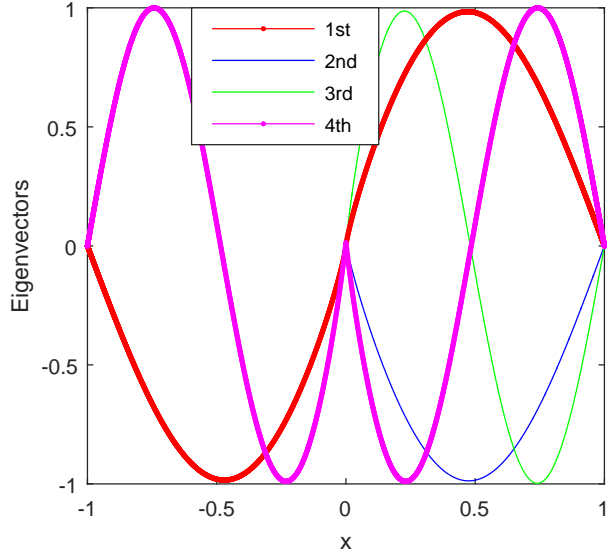}
\caption{\fontsize{9}{9}\selectfont The first four eigenvectors of problem (\ref{B_HWV}) computed by Chebfun when $\tau:=0$ and $\nu:=1/3$.}
\label{Bessel_0_1_3}
\end{figure}
For $\nu \in \left[1/3, 1/2\right]$ in \cite{HWV} it is shown that $\lambda _{n}\left( \tau \right) $, $n\in 
\mathbb{N}$ are analytic in the
domain  $
\mathbb{C}\diagdown \left[ -1,1\right] $ and at $\tau :=\infty $ and these functions
are also analytic on the segments $\tau \pm i0$ for $-1<\tau \,<1$ and they
can be extended continuously into $\tau =\pm 1$ respectively.

With the following simple Chebfun code
\begin{verbatim}
    x=chebfun('x',[-1, 1],'splitting','on'); nu=1/3;tau=0; %tau=1;
L=chebop(-1,1); 
L.op =@(x,y)-((x-tau)^2)*diff(y,2)+((nu^2)-1/4)*y; 
L.rbc=0;  L.lbc=0;
M=chebop(-1,1); M.op=@(x,y) ((x-tau)^2)*y; M.bc=[ ];
[V,D]=eigs(L,M,40); sort(diag(D))
\end{verbatim}
we have computed the eigenpairs of problem (\ref{B_HWV}) for various $\tau$ and $\nu$. They are reported in the first two columns of Table \ref{cross_eigs}. 

It is very important and honest to note that for the values of $\tau$ close to zero Chebfun displays the following warning: \textit{Maximum dimension reached. Solution may not have converged.} This warning prompted us to look for a way to validate these results. Thus we have solved the same problem with the classic ChC. For an order of approximation $N:=2028$, i.e., half of that used by Chebfun, that method produced the third column of the Table \ref{cross_eigs}. A simple comparison of the first and the third column of the table leads us to believe that the results are correct.

\begin{table} \centering%
\begin{tabular}{|l|l|l||l|}
\hline
$j\:$ &\:$\lambda_{j}$ for $\tau:=0$ Chebfun & \:$\lambda_{j}$ for $\tau:=1/4$ Chebfun &\:$\lambda_{j}$ for $\tau:=0$ ChC \\ \hline
$0\:$ &$\:8.427067009456547e+00 $& $\:5.400563866142070e+00 $ & $\: 8.427945713285165e+00$\\ \hline
$1\:$ &$\:8.571978456544725e+00 $& $\:1.503383426547553e+01 $ & $\: 8.661777916480734e+00$\\ \hline
$2\:$ &$\:3.640104840502095e+01 $& $\:2.332141941007685e+01 $ & $\: 3.640525136075537e+01$\\ \hline
$3\:$ &$\:3.689758359379723e+01 $& $\:5.388134332882591e+01 $ & $\: 3.720851434888583e+01$\\ \hline
$4\:$ &$\:8.411229915649267e+01 $& $\:6.489190302110852e+01 $ & $\: 8.412057813169606e+01$\\ \hline
$5\:$ &$\:8.511762862999629e+01 $& $\:9.708229958053987e+01 $ & $\: 8.575213457278441e+01$\\ \hline
$6\:$ &$\:1.515639240123626e+02 $& $\:1.498879163534515e+02 $ & $\: 1.515775098405456e+02$\\ \hline
$7\:$ &$\:1.532164600154869e+02 $& $\:1.529282181428597e+02 $ & $\: 1.542668794303226e+02$\\ \hline
\hline
\end{tabular}%
\caption{The first eight eigenvalues of eigenproblem (\ref{B_HWV}) computed by Chebfun for $\tau:=0$ and $\tau:=1/4$
(first two columns). In the third column we report the eigenvalues of the same problem computed by ChC when $\tau:=0$. In both methods we have assumed $\nu=1/3$. }
 \label{cross_eigs}%
\end{table}
The first four eigenvectors of problem (\ref{B_HWV}) computed by Chebfun when $\tau:=0$ and $\nu:=1/3$ are displayed in Fig. \ref{Bessel_0_1_3}. From this figure it is clear that, in pairs, i.e., for
1 and 2, 3 and 4, etc. the eigenvectors are identical over the range $\left[-1, 0\right]$. More important and at the same time very interesting, they become perfectly symmetrical with respect to the horizontal axis over the interval $\left[0, 1\right]$.
\begin{figure}
\centering
\includegraphics[scale=0.85]{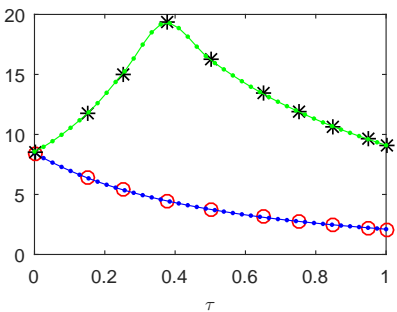}
\caption{\fontsize{9}{9}\selectfont The dependence of the first two sets of eigenvalues $\lambda_{0}$ (lower curve) and $\lambda_{1}$ of problem (\ref{B_HWV}) on $\tau$ when $\nu:=1/3$.}
\label{eig_tau_dep}
\end{figure}

The dependence of the first two sets of eigenvalues $\lambda_{0}$ and $\lambda_{1}$ of problem (\ref{B_HWV}) on $\tau$ is depicted in Fig. \ref{eig_tau_dep}. The lower curve corresponds to $\lambda_{0}\left(\tau\right)$, $\tau \in \left[0, 1\right].$ This figure has to be seen with a symmetric part on $\left[-1, 0\right]$. In order to smooth these curves we have used the MATLAB \texttt{pchip} code that performs a piecewise Hermite cubic interpolation. 

In our opinion this figure suggests that in origin occurs the so-called \textit{phenomenon of level crossing of eigenvalues} as described in \cite{BO} p. 350 (see also \cite{LNTex}). Due to the scale at which the curves are drawn the first two eigenvalues $\lambda_{0}$ and $\lambda_{1}$, from the first column of the Table \ref{cross_eigs}, appear to be superimposed but actually they are only very close.

Along with the results displayed in Table \ref{cross_eigs} we can answer to the question formulated in the title of paper \cite{Teytel}. Thus, well, close to the origin multiple eigenvalues are not rare at all for the Bessel problem (\ref{B_HWV}). Moreover,this result is in the spirit of Theorem 6.6 from \cite{Teytel} on degeneracies in the spectra of self-adjoint operators (of Schr\"odinger type).

\section{Concluding remarks}\label{conclusions}

Regarding the regular problems, Chebfun is unbeatable in terms of accuracy, computation speed, and the information they provide on the accuracy of computational process. It displays the optimal approximation order of unknowns (eigenvectors) and how and to what extent their Chebyshev coefficients decrease. It also specifies the degree to which the boundary conditions are satisfied. 

As for the singularly perturbed problems, the situation is not so offering. In this case, Chebfun also leaves a bit of room for the usual (classical) spectral methods. However, approaching these problems in parallel, using Chebfun as well as the classical spectral methods, gives greater confidence in the accuracy of their results. Regarding ChC method, using the so-called relative drift  to some parameters we ensure the numerical stability of the numerical process.

The first open problem that comes to mind, and for which there are some results in progress, is the simultaneous dependence of the eigenvalues of Bessel equation on both parameters $\tau$ and $\nu$. This could shed more light on the singularity in origin.


\begin{thebibliography}{99}

\bibitem{BEZ} P. B. Bailey, W. N. Everitt and A. Zettl, \emph{Computing Eigenvalues of Singular Sturm-Liouville Problems}, Results Math.,20 (1991) pp. 391--423.

\bibitem{BGKZ} P. Bailey, B. Garbow, H. Kaper, and A. Zettl, \emph{Algorithm 700: A FORTRAN software
package for Sturm-Liouville problems}, ACM T. Math. Software, 17 (1991), pp. 500--501.

\bibitem{BO} C. Bender and S. Orszag, \emph{Advanced Mathematical Methods for Scientists and Engineers}, McGraw-Hill, New York 1978.

\bibitem{Boyd} J. P. Boyd, \emph{Traps and Snares in Eigenvalue Calculations with Application to Pseudospectral 
Computations of Ocean Tides in a Basin Bounded by Meridians}, J. Comput. Phys. 126 (1996) pp. 11--20.

\bibitem{Castro}  H. Castro, and H. Wang, \emph{A singular Sturm-Liouville equation under
homogeneous boundary conditions}, J. Func. Anal. 261 (2011) pp. 1542--1590.

\bibitem{DBT} T. A. Driscoll, F. Bornemann, and L. N. Trefethen,
\emph{The CHEBOP System for Automatic Solution of Differential Equations}, BIT 48 (2008) pp. 701-–723.

\bibitem{DHT} T. A. Driscoll, N. Hale, and L. N. Trefethen (eds.), \emph{Chebfun Guide}, Pafnuty Publications, Oxford, 2014

\bibitem{DHT1} T. A. Driscoll, N. Hale, and L. N. Trefethen, \emph{Chebfun-numerical computing with functions}. Software available at http://www.chebfun.org.

\bibitem{EGZ} W.N. Everitt, J. Gunson, and A. Zettl, \emph{Some comments on Sturm-Liouville eigenvalue problems
with interior singularities}, Journal of Applied Mathematics and Physics (ZAMP) 38 (1987) pp. 813--838.

\bibitem{GigIsp} C. I. Gheorghiu and I. S. Pop, \emph{A Modified Chebyshev-Tau Method for a Hydrodynamic Stability Problem}, Proceedings of the International
Conference on Approximation and Optimization (Romania) - ICAOR
Cluj-Napoca, July 29 - August 1, 1996. Volume II, pp. 119--126.

\bibitem{CIGHPR}  C.I. Gheorghiu, M.E. Hochstenbach,  B. PLestenjak, and J. Rommes, \emph{Spectral collocation solutions to multiparameter Mathieu's system}, Appl. Math. Comput. 218 (2012) pp. 11990–12000.

\bibitem{Cig14}
C.~I. Gheorghiu, \emph{Spectral Methods for Non-Standard Eigenvalue Problems. Fluid and Structural Mechanics and Beyond}, Springer-Verlag, Cham Heidelberg New-York Dondrecht London, 2014.
\bibitem{Cig17} C.~I. Gheorghiu,\emph{On the numerical treatment of the eigenparameter dependent boundary conditions}, Numer.
Algor. 77 (2018) pp. 77–-93.

\bibitem{Cig18} C. I. Gheorghiu, \emph{Spectral Collocation Solutions to Problems on Unbounded Domains}, Casa C\u{a}r\c{t}ii de \c{S}tiin\c{t}\u{a} Publishing House, Cluj-Napoca, 2018

\bibitem{JH} J. Hoepffner, \emph{Implementation of boundary
conditions}. Software available at http://www.lmm.jussieu.fr/$\sim$
hoepffner/boundarycondition.pdf

\bibitem{MP} M. Marletta, and J. D. Pryce, \emph{LCNO Sturm-Liouville problems—Computational
difficulties and examples}, Numer. Math. 69  (1995) pp. 303–-320.

\bibitem{OT} S. Olver and A. Townsend, \emph{A fast and well-conditioned spectral method}, SIAM Rev., 55 (2013), pp. 462–-489.

\bibitem{Pleijel} A. Pleijel, \emph{On the  Boundary Conditions for the Legendre Polynomials}, Ann. Acad. Sci. Fenn. A1 2 (1976) pp. 397--408.

\bibitem{PGH} B. Plestenjak, C. I. Gheorghiu, and M. E.Hochstenbach, \emph{Spectral collocation for multiparameter eigenvalue problems arising from separable boundary value problems}, J. Comput. Phys. 298 (2015) pp. 585–-601.

\bibitem{PF}  S. Pruess, and C. T. Fulton, \emph{Mathematical Software for
Sturm-Liouville Problems}, ACM T. Math. Software 19
(1993) pp. 360--376. 

\bibitem{PFX} S. Pruess, C. T. Fulton, and Y. Xie, \emph{An Asymptotic Numerical Method for a Class of Singular Sturm-Liouville Problems}, SIAM J. Numer. Anal. 32 (1995) pp. 1658--1676.

\bibitem{PM} J. Pryce and M. Marletta, \emph{A new multi-purpose software package for Schrodinger and
Sturm--Liouville computations}, Comput. Phys. Comm., 62 (1991), pp. 42--54.

\bibitem{Pr} J. D. Pryce, \emph{A Test Package for Sturm-Liouville
Solvers}, ACM T. Math. Software 25 (1999) pp. 21--57. 

\bibitem{BDS} B. D. Shizgal, \emph{Pseudospectral Solution of the Fokker-Planck Equation with Equilibrium Bistable States: the Eigenvalue spectrum and the Approach to Equilibrium}, J. Stat. Phys. (2016)

\bibitem{S}  C.A. Stuart, \emph{On the spectral theory of a tapered rod},  Proc. Roy. Soc. Edinburgh Sect. A
132 (2002) pp. 729--764.

\bibitem{Teytel} M. Teytel, \emph{How Rare Are Multiple Eigenvalues?}, Comm. Pure Appl. Math. 52 (1999) pp. 917--934.

\bibitem{LNTex} L. N. Trefethen, \emph{Analyticity at eigenvalue near-crossings}, Software available at https://www.chebfun.org/examples/linalg/CrossingsAnalyticity.html.

\bibitem{TBD} L.~N. Trefethen, A. Birkisson, and T. A. Driscoll, \emph{Exploring ODEs}, SIAM, Philadelphia, 2018.

\bibitem{LNT} L. N. Trefethen, \emph{Approximation Theory and Approximation Practice, Extended Edition},
SIAM, Philadelphia, 2019

\bibitem{HWV} H. W. Volkmer, \emph{Eigenvalue problems for Bessel’s equation and zero-pairs of Bessel functions}, 
Stud. Sci. Math. Hung. 35 (1999) pp.261--280.

\bibitem{WR} J. A. C. Weideman, and S. C. Reddy, \emph{A MATLAB
Differentiation Matrix Suite}, ACM T. Math. Software,
26, (2000) pp. 465--519. 

\bibitem{GvW} G. von Winckel, \emph{Fast Chebyshev Transform (1D)}, Software available at https://www.mathworks.com/matlabcentral/fileexchange/4591-fast-chebyshev-transform-1d.
\end{thebibliography}
\end{document}